\documentclass[12pt,reqno]{article}

\usepackage[usenames]{color}
\usepackage{amssymb}
\usepackage{amsmath}
\usepackage{amsthm}
\usepackage{amsfonts}
\usepackage{amscd}
\usepackage{graphicx}
\usepackage[utf8]{inputenc}
\usepackage[T1]{fontenc}

\usepackage[colorlinks=true,
linkcolor=webgreen,
filecolor=webbrown,
citecolor=webgreen]{hyperref}

\definecolor{webgreen}{rgb}{0,.5,0}
\definecolor{webbrown}{rgb}{.6,0,0}

\usepackage{color}

\usepackage{graphics}
\usepackage{latexsym}

\begin{document}

\theoremstyle{plain}
\newtheorem{theorem}{Theorem}
\newtheorem{corollary}[theorem]{Corollary}
\newtheorem{lemma}[theorem]{Lemma}
\newtheorem{proposition}[theorem]{Proposition}

\theoremstyle{definition}
\newtheorem{definition}[theorem]{Definition}
\newtheorem{example}[theorem]{Example}
\newtheorem{conjecture}[theorem]{Conjecture}

\theoremstyle{remark}
\newtheorem{remark}[theorem]{Remark}

\begin{center}
\vskip 1cm{\LARGE\bf Is This a New Class of Matrices?
\vskip 1cm}
\large
Jovan Miki\'{c}\\
University of Banja Luka\\
Faculty of Technology\\
Bosnia and Herzegovina\\
\href{mailto:jovan.mikic@tf.unibl.org}{\tt jovan.mikic@tf.unibl.org} \\
\end{center}

\vskip .2in

\begin{abstract}

We consider a new class of matrices associated to a real square matrix $A$  and to a vector $\vec{c} \in \{-1,1\}^n$ such that $c_1=1$ by using a map  $\varphi_{\vec{c}}$ which turns out to be a conjugation of a matrix $A$ by a signature matrix. It is shown that every such matrix is similar and congruent to a matrix $A$ and that they have same permanental polynomials.  There are $2^{n-1}$ maps $\varphi_{\vec{c}}$ and they form an abelian group under the composition of maps isomorphic to the group $({\mathbb{Z}_2}^{n-1}, +)$. A decomposition of matrices, on  a symmetric and antisymmetric matrix under a map $\varphi_{\vec{c}}$, is considered. Particularly, it is shown that sum of all principal minors of the order two of a matrix $A$ is equal to the sum of all principal minors of the order two of their symmetric and antisymmetric parts.  It is shown that any symmetric matrix  and any antisymmetric matrix under the map $\varphi_{\vec{c}}$ are simultaneously permutation similar to certain block matrices which have two blocks.  Finally, for a fixed matrix $A$, it is proved that the number of  different matrices $\varphi_{\vec{c}}(A)$ is $2^{n-t}$, where $t$ is the number of connected components of the graph $G$ whose adjacency matrix  is $A$.

\end{abstract}


\noindent \textbf{2020} {\it \textbf{Mathematics Subject Classification}}:  15B35, 15B57

\section{Introduction}\label{sec:1}

 Let $A$  be a square real matrix, and let  $P$ be a permutation matrix of the same order as matrix  $A$. It is known \cite[Eq.~(2), p.\ 274]{RM} that a matrix $P^{-1}\cdot A \cdot P$ is similar and congruent to a matrix $A$, and they have same permanental polynomials. It is well-known \cite[Chapter 2]{SP} that any permutation matrix must be an orthogonal matrix:  $P^{-1}=P^{T}$.

Also, it is known that a transpose matrix  $A^T$ is similar to a matrix $A$ and  they have same permanental polynomial. Not so long ago, it is proved  \cite{DjokoIkramov, RHVS} that  $A^T$ is congruent to $A$.

In this paper, we introduce a new class of matrices (up to our knowledge)  which are similar and congruent to a matrix $A$, and they have same permanental polynomials as a matrix $A$.

 Let  $A=\{a_{ij}:1\leq i,j \leq n\}\in M_{n,n}(\mathbb{R})$ be an arbitrary square real matrix of an order $n$,
 and let $\vec{c}=(c_i)_{i=1}^{n} \in \{-1,1\}^n$  be a vector such that $c_1=1$.

Let us consider a map $\varphi_{\vec{c}}:M_{n,n}(\mathbb{R}) \to M_{n,n}(\mathbb{R})$, as follows:
\begin{equation}\label{eq:1}
(\varphi_{\vec{c}}(A))_{ij}=c_i a_{ij} c_j\text{.}
\end{equation}

Let  $tr(A)$,  $|A|$, $perm(A)$  denote the trace, the determinant, and  the permanent of a matrix $A$, respectively. 
We begin with the following two lemmas:
\begin{lemma}\label{l:1}
The following relations are true:
\begin{align}
tr(\varphi_{\vec{c}}(A))&=tr(A)\text{,}\label{eq:2}\\
|\varphi_{\vec{c}}(A)|&=|A|\text{,}\label{eq:3}\\
perm(\varphi_{\vec{c}}(A))&=perm(A)\text{.}\label{eq:4}
\end{align}
\end{lemma}

Let $A(i_1,i_2,\ldots,i_k)$ denote the principal minor of a matrix $A$ of an order $k$, where $1 \leq i_1 <i_2< \ldots < i_k \leq n$ and $1\leq k \leq n$. Similarly,  let $A[i_1,i_2,\ldots,i_k]$ denote the principal permanent of a matrix $A$ of an order $k$. Furthermore,  let $p_A(\lambda)=|A-\lambda I|$ denote the characteristic polynomial  of a matrix $A$, and let $q_A(\lambda)=per(A-\lambda I)$ denote the permanental polynomial of a matrix $A$.

\begin{lemma}\label{l:2}

The following relations are  true:
\begin{align}
\varphi_{\vec{c}}(A)(i_1,i_2,\ldots,i_k)&=A(i_1,i_2,\ldots,i_k)\text{,}\label{eq:5}\\
p_{\varphi_{\vec{c}(A)}}(\lambda)&=p_A(\lambda)\label{eq:6}\\
\varphi_{\vec{c}}(A)[i_1,i_2,\ldots,i_k]&=A[i_1,i_2,\ldots,i_k]\label{eq:7.1}\\
q_{\varphi_{\vec{c}(A)}}(\lambda)&=q_A(\lambda)\label{eq:7.2}\text{.}
\end{align}
\end{lemma}

Also, the following corollary is true:

\begin{corollary}\label{cor:1}
\begin{equation}
rank(\varphi_{\vec{c}}(A))=rank(A)\text{.}\label{eq:7}
\end{equation}
\end{corollary}

We assert that:

\begin{theorem}\label{t:1}
Matrix $\varphi_{\vec{c}}(A)$ is similar and congruent to a matrix $A$, and they have the same permanental polynomial.
\end{theorem}

\begin{theorem}\label{t:2}
There are, in total, $2^{n-1}$ different maps $\varphi_{\vec{c}}$.
The set $\varPsi_n$ of all maps $\varphi_{\vec{c}}$ form an abelian group under the  composition of maps. This group is isomorphic with the group $({\mathbb{Z}_2}^{n-1}, +)$.
\end{theorem}

Let $Sym$  and $AntiSym$ denote sets of  symmetric and antisymmetric matrices of an order $n$ over the field of real numbers, respectively. Note that we use the word antisymmetric instead of skew-symmetric.

It is well-known that:

\begin{equation}\notag
Sym \oplus AntiSym=M_{n,n}(\mathbb{R})\text{,}
\end{equation}
where $\oplus$ denotes the direct sum.

We consider the following decomposition of matrices:
\begin{definition}\label{def:1}
We say that matrix $A \in M_{n,n}(\mathbb{R})$ is symmetric under the map $\varphi_{\vec{c}}$ if  $\varphi_{\vec{c}}(A)=A$.  Matrix $A$ is antisymmetric  under the map $\varphi_{\vec{c}}$ if  $\varphi_{\vec{c}}(A)=-A$.
Let $Sym_{\vec{c}}$ denote the set of all symmetric matrix under the map $\varphi_{\vec{c}}$,  and let $AntiSym_{\vec{c}}$ denote the set of all antisymmetric matrix under the map $\varphi_{\vec{c}}$.
\end{definition}

We assert that:

\begin{theorem}\label{t:3}
The algebraic structures $(Sym_{\vec{c}},+)$ and $(AntiSym_{\vec{c}},+)$ form vector subspaces of the vector space $(M_{n,n}(\mathbb{R}),+)$ over the field of real numbers; where $+$ denotes addition of matrices.

Let $A\in M_{n,n}(\mathbb{R})$. Let  us define $A_S(\vec{c})$ and $A_{AS}(\vec{c})$  matrices, as follows:
\begin{equation}\label{eq:8}
(A_S(\vec{c}))_{ij}=
\begin{cases}
a_{ij}, &\text{ if $c_i\cdot c_j=1$ }\\
0&\text{ if $c_i\cdot c_j=-1$}
\end{cases}
\end{equation}

\begin{equation}\label{eq:9}
(A_{AS}(\vec{c}))_{ij}=
\begin{cases}
0, &\text{ if $c_i\cdot c_j=1$ }\\
a_{ij}&\text{ if $c_i\cdot c_j=-1$ }
\end{cases}
\end{equation}

Then $A_S(\vec{c})\in Sym_{\vec{c}}$, $A_{AS}(\vec{c})\in AntiSym_{\vec{c}}$, and $A=A_S(\vec{c})+A_{AS}(\vec{c})$.

We obtain that:
\begin{equation}\label{eq:10}
Sym_{\vec{c}} \oplus AntiSym_{\vec{c}}=M_{n,n}(\mathbb{R})\text{.}
\end{equation}
\end{theorem}

By Theorem \ref{t:3}, directly follows that $(A_S(\vec{c}))_{ii}=a_{ii}$ and $A_{AS}(\vec{c})_{ii}=0$, for all $1\leq i \leq n$. 
\begin{lemma}\label{l:3}
The algebraic structure $(Sym_{\vec{c}},+,\cdot)$ is a subring of a ring $(M_{n,n}(\mathbb{R}),+,\cdot)$, where $+$ and $\cdot$ denote addition and standard multiplication of matrices, respectively.
\end{lemma}

Note that this decomposition of matrices has the following property:

\begin{theorem}\label{t:4}
Let $A\in M_{n,n}(\mathbb{R})$. Then the sum of all principal minors of the order two of matrices $A_S(\vec{c})$ and $A_{AS}(\vec{c})$ is equal to  the sum of all principal minors of the order two of a matrix $A$. Similarly, the sum of all principal permanents of the order two of matrices $A_S(\vec{c})$ and $A_{AS}(\vec{c})$ is equal to  the sum of all principal permanent of the order two of a matrix $A$.
\end{theorem}

Let $[\lambda^{k}]p(\lambda)$ denote the coefficient of a polynomial $p$  by the $\lambda^k$.

Theorem \ref{t:4} can be reformulated as:
\begin{align}\
[\lambda^{n-2}]p_A(\lambda)&=[\lambda^{n-2}]p_{(A_S(\vec{c}))}(\lambda)+[\lambda^{n-2}]p_{(A_{AS}(\vec{c}))}(\lambda)\label{eq:11}\\
[\lambda^{n-2}]q_A(\lambda)&=[\lambda^{n-2}]q_{(A_S(\vec{c}))}(\lambda)+[\lambda^{n-2}]q_{(A_{AS}(\vec{c}))}(\lambda)\text{.}\label{eq:11.1}
\end{align}
The same property, described in Theorem \ref{t:4}, is true also for a  decomposition of matrices on symmetric and antisymmetric matrices.

\begin{theorem}\label{t:5}
Let $A\in M_{n,n}(\mathbb{R})$, and let $A^{T}$ denote the transpose matrix of $A$.  The sum of all principal minors of the order two of matrices $A_S=\frac{1}{2}(A+A^{T})$ and $A_{AS}=\frac{1}{2}(A-A^{T})$ is equal to  the sum of all principal minors of the order two of a matrix $A$. Similarly, the sum of all principal permanents of an order two of matrices $A_S$ and $A_{AS}$ is equal to the sum of all principal permanents of the order two of a matrix $A$.
\end{theorem}

Theorem \ref{t:5} can be reformulated as:
\begin{align}
[\lambda^{n-2}]p_A(\lambda)&=[\lambda^{n-2}]p_{(A_S)}(\lambda)+[\lambda^{n-2}]p_{(A_{AS})}(\lambda)\label{eq:12}\\
[\lambda^{n-2}]q_A(\lambda)&=[\lambda^{n-2}]q_{(A_S)}(\lambda)+[\lambda^{n-2}]q_{(A_{AS})}(\lambda)\label{eq:12.1}\text{.}
\end{align}

\begin{definition}\label{def:2.1}
Let $A$ and  $B \in M_{n,n}(\mathbb{R})$. Matrix $B$ is permutation similar to matrix $A$ if  $B=P^{-1}\cdot A\cdot P$, where $P$ is a permutation matrix.
\end{definition}

We assert that any symmetric matrix under the map $\varphi_{\vec{c}}$ is  permutation similar to a block-diagonal matrix which has two blocks.
\begin{definition}\label{def:2}

Let $\vec{c}\in \{-1,1\}^n$ such that $c_1=1$, and let $A \in  Sym_{\vec{c}}$. Let us assume that $\vec{c}$ has exactly $r$-coordinates which are equal to one, where $1\leq r \leq n$.
Let $I_r$ denote the set $\{i_1,i_2,\ldots, i_r\}$ such that $c_{i_k}=1$ for all $1\leq k \leq r$; where $1\leq i_1<i_2<\ldots < i_r\leq n$.
Let $J_{n-r}$ denote the set $\{j_1,j_2,\ldots, j_{n-r}\}$ such that $c_{j_k}=-1$ for all $1\leq k \leq n-r$;  where $1< j_1<j_2<\ldots < j_{n-r}\leq n$.

Let us define a matrix $D=\{d_{l m}\} \in M_{r,r}(\mathbb{R})$ such that
$d_{l m}=a_{i_l i_m}$ for all $1\leq l,m \leq r$.

Furthermore, let us define a matrix $E=\{e_{l m}\} \in M_{n-r,n-r}(\mathbb{R})$ such that
$e_{l m}=a_{j_l j_m}$ for all $1\leq l,m \leq n-r$.

\end{definition}

\begin{theorem}\label{t:6}
 Let $\vec{c}\in \{-1,1\}^n$ such that $c_1=1$, and let $A \in  Sym_{\vec{c}}$. Let $D$ and $E$ be matrices from Definition \ref{def:2}. Then matrix $A$ is  permutation similar to the block-diagonal  matrix $diag(D,E)$.
\end{theorem}

\begin{corollary}\label{c:4}
 Let $\vec{c}\in \{-1,1\}^n$ such that $c_1=1$, and let $A \in  Sym_{\vec{c}}$. Let $D$ and $E$ be matrices from Definition \ref{def:2}. Then we have:
\begin{align}
p_A(\lambda)&=p_D(\lambda)\cdot p_E(\lambda)\label{eq:16.1}\\
|A|&=|D|\cdot|E|\label{eq:16.2}\\
perm(A)&=perm(D)\cdot perm(E)\text{.}\label{eq:16.3}
\end{align}
\end{corollary}

Also, we assert that any antisymmetric matrix under the map $\varphi_{\vec{c}}$ is permutation similar to an another block matrix which has two blocks.

\begin{definition}\label{def:3}

Let $\vec{c}\in \{-1,1\}^n$ such that $c_1=1$, and let $A \in AntiSym_{\vec{c}}$. 
Let $I_r$ and $J_{n-r}$ be sets from Definition \ref{def:2}.

Let us define a matrix $F=\{f_{l m} \}\in M_{r,n-r}(\mathbb{R})$ such that
$f_{l m}=a_{i_l j_m}$ for all $1\leq l \leq r$ and $1\leq m\leq n-r$.

Let us define a matrix $G=\{g_{l m}\} \in M_{n-r,r}(\mathbb{R})$ such that
$g_{l m}=a_{j_l i_m}$ for all $1\leq l \leq n-r$ and $1\leq m\leq r$.

Let $H$ denote the following block matrix:

\begin{equation}\label{eq:16.4}
\begin{pmatrix}
O_{r}&F\\
G&O_{n-r}
\end{pmatrix}
\end{equation}
; where $O_r\in M_{r,r}(\mathbb{R})$ and $O_{n-r}\in M_{n-r,n-r}(\mathbb{R})$ are zero matrices.
\end{definition}

\begin{theorem}\label{t:7}
 Let $\vec{c}\in \{-1,1\}^n$ such that $c_1=1$, and let $A \in  AntiSym_{\vec{c}}$. Let $H$ be a block matrix from the Eq.~ \textup{(}\ref{eq:16.4}\textup{)}. Then matrix $A$ is permutation similar to $H$.
\end{theorem}

\begin{corollary}\label{c:7}
 Let $\vec{c}\in \{-1,1\}^n$ such that $c_1=1$, and let $A \in  AntiSym_{\vec{c}}$. Let us assume that $\vec{c}$ has exactly $r$-coordinates which are equal to one, where $1\leq r \leq n$.
Then $|A|=0$ and $perm(A)=0$ except for the case where $n$ is even and $r=\frac{n}{2}$. In that case, the following two equations are true:
\begin{align}
|A|&=(-1)^n\cdot |F|\cdot |G|\label{eq:16.5}\\
perm(A)&=perm(F)\cdot perm (G)\label{eq:16.6}
\end{align}
\end{corollary}

Finally, let us assume that a matrix $A \in M_{n,n}(\mathbb{R})$ is fixed. Let $\vec{c}\in \{-1,1\}^n$ such that $c_1=1$. A question that arises is: ``How many  are different matrices $\varphi_{\vec{c}}(A)$ ?''.

\begin{definition}\label{def:4}
Let $G$ be a simple undirected graph with $n$ vertices such that its adjacency matrix is the matrix A, i.e., vertices $v_i$ and $v_j$ of a graph $G$  are adjacent if, at least, one of the two elements  $a_{ij}$ or $a_{ji}$ of the matrix $A$ is non-zero.
\end{definition}
\begin{theorem}\label{t:8}
Let $t$ be a number of connected components of the graph $G$ from Definition \ref{def:4}.
Then there are $2^{n-t}$ different matrices $\varphi_{\vec{c}}(A)$.
\end{theorem}

\begin{corollary}\label{c:8}
Let  a matrix $A \in M_{n,n}(\mathbb{R})$ be fixed, and let $\vec{c}\in \{-1,1\}^n$ such that $c_1=1$.
There are  $2^{t-1}$ matrices such that $\varphi_{\vec{c}}(A)=A$.
\end{corollary}

\section{A Proof of Lemma \ref{l:1}}\label{sec:2} 

\subsection{A Proof of the Eq.~(\ref{eq:2})}\label{sec:2.1}

Let $i$ be fixed natural number such that $1\leq i \leq n$.
By the Eq.~(\ref{eq:1}) and the fact that  $\vec{c} \in \{-1,1\}^n$, it follows that
\begin{equation}\label{eq:13}
(\varphi_{\vec{c}}(A))_{ii}={c_i}^2 a_{ii}=a_{ii}\text{.}
\end{equation}

By using the Eq.~(\ref{eq:13}), it follows that
\begin{align*}
tr(\varphi_{\vec{c}}(A))&=\sum_{i=1}^{n}(\varphi_{\vec{c}}(A))_{ii}\\
&=\sum_{i=1}^{n}a_{ii}\\
&=tr(A)\text{.}
\end{align*}
The last equation above proves the Eq.~(\ref{eq:2}).

\subsection{Proofs of the Eqns.~(\ref{eq:3}) and (\ref{eq:4})}\label{sec:2.2}
Let $i$  and $j$ be natural numbers such that $1\leq i ,j\leq n$.
By the Eq.~(\ref{eq:1}), it follows that  every $i$-th row of a matrix $\varphi_{\vec{c}}(A)$ is a multiple of an integer $c_i$. Note that $c_i$ is on the left-side from $a_{ij}$.
Similarly, every $j$-th column of a matrix $\varphi_{\vec{c}}(A)$ is a multiple of an integer $c_j$. Note that $c_j$ is on the right-side from $a_{ij}$.

From every  $i$-th row of a $|\varphi_{\vec{c}}(A)|$, we can extract an integer $c_i$ in front of the determinant. Similarly, from every  $j$-th column of a $|\varphi_{\vec{c}}(A)|$ we can extract an integer $c_j$ in front of the determinant.

We have gradually:

\begin{align*}
|\varphi_{\vec{c}}(A)|&=\bigl{(}\prod_{i=1}^{n}c_i\bigr{)}\cdot \bigl{(}\prod_{j=1}^{n}c_j\bigr{)}\cdot|A|\\
&=\bigl{(}\prod_{i=1}^{n}c_i\bigr{)}^2\cdot |A|\\
&=\bigl{(}\prod_{i=1}^{n}{c_i}^2\bigr{)}\cdot |A| \\
&=|A|\text{.} && ( \text{by using } \vec{c} \in \{-1,1\}^n )
\end{align*}

This proves the Eq.~(\ref{eq:3}).

The proof of the Eq.~(\ref{eq:4}) is similar to the proof of the Eq.~(\ref{eq:3}). Just replace $||$ with $perm()$ in the previous proof.  This  proof is correct due to the fact that multiplying any single row or column of a square matrix $A$ by a real number $s$ changes $perm(A)$ to $s\cdot perm(A)$. 
This completes the proof of Lemma \ref{l:1}.

\section{A Proof of Lemma \ref{l:2}}\label{sec:3}
\subsection{A Proof of the Eq.~(\ref{eq:5})}\label{sec:3.1}

Let $l$ and $m$ be natural numbers such that $1\leq l,m \leq k$.

The $l$-th row of  $\varphi_{\vec{c}}(A)(i_1,i_2,\ldots,i_k)$ is a multiple of an integer $c_{i_l}$.
Similarly, $m$-th column of  $\varphi_{\vec{c}}(A)(i_1,i_2,\ldots,i_k)$ is a multiple of an integer $c_{i_m}$.

We have gradually:

\begin{align*}
\varphi_{\vec{c}}(A)(i_1,i_2,\ldots,i_k)&=\bigl{(}\prod_{l=1}^{k}c_{i_l}\bigr{)}\cdot \bigl{(}\prod_{m=1}^{k}c_{i_m}\bigr{)}\cdot A(i_1,i_2,\ldots,i_k)\\
&=\bigl{(}\prod_{l=1}^{k}c_{i_l}\bigr{)}^2\cdot A(i_1,i_2,\ldots,i_k)\\
&=\bigl{(}\prod_{l=1}^{k}{c_{i_l}}^2\bigr{)}\cdot  A(i_1,i_2,\ldots,i_k) \\
&=  A(i_1,i_2,\ldots,i_k)\text{.} && ( \text{by using } \vec{c} \in \{-1,1\}^n )
\end{align*}

\subsection{A Proof of the Eq.~(\ref{eq:6})}\label{sec:3.2}
Let $k$ be a natural number such that $k\leq n$.

It is well-known \cite[p.\ 117]{MJ} that:

\begin{equation}\label{eq:14}
[\lambda^{n-k}] p_A(\lambda)=(-1)^{n-k}\underset{1\leq i_1<i_2<\ldots <i_k\leq n}\sum A(i_1,i_2,\ldots,i_k)\text{;}
\end{equation}
where summation is over all principal minors of a matrix $A$ of an order $k$.

The proof of the Eq.~(\ref{eq:6}) follows by Eqns.~(\ref{eq:5}) and (\ref{eq:14}).

This completes the proof of Lemma \ref{l:2}.

\subsection{Proofs of Eqns.~(\ref{eq:7.1}) and (\ref{eq:7.2})}

A proof of Eq.~(\ref{eq:7.1}) is analogue with the proof of Eq.~(\ref{eq:5}). Just replace principal minors $\varphi_{\vec{c}}(A)(i_1,i_2,\ldots,i_k)$  and $A(i_1,i_2,\ldots,i_k)$  with principal permanents $\varphi_{\vec{c}}(A)[i_1,i_2,\ldots,i_k]$  and $A[i_1,i_2,\ldots,i_k]$, respectively.  This proof is correct due to the fact that multiplying any single row or column of a square matrix $A$ by a real number $s$ changes $perm(A)$ to $s\cdot perm(A)$.

Furthermore, it can be shown that \cite[Eq.~(1), p.\ 273]{RM}:

\begin{equation}\label{eq:15.1}
[\lambda^{n-k}] q_A(\lambda)=(-1)^{n-k}\underset{1\leq i_1<i_2<\ldots <i_k\leq n}\sum A[i_1,i_2,\ldots,i_k]\text{;}
\end{equation}
where summation is over all principal permanents of a matrix $A$ of an order $k$.

The proof of the Eq.~(\ref{eq:15.1}) is analogue to the proof of the Eq.~(\ref{eq:14}).

Therefore, the proof of the Eq.~(\ref{eq:7.2}) follows from Eqns.~(\ref{eq:7.1}) and (\ref{eq:15.1}).

\section{A Proof of Corollary \ref{cor:1}}\label{sec:4}

The matrix $\varphi_{\vec{c}}(A)$ can be obtained from the matrix $A$ by multiplication of $i$-th row by non-zero integer $c_i\in \{-1,1\}$,  and by  multiplication of $j$-th column with non-zero integers $c_j\in \{-1,1\}$; where $1\leq i ,j\leq n$. There are, in total, $2n$ multiplications.
Since multiplication of row (column) by non-zero real number does not change the rank of a matrix, the Eq.~(\ref{eq:7}) follows.
This completes the proof of Corollary \ref{cor:1}.

\section{A Proof of Theorem \ref{t:1}}\label{sec:5}

Let $\vec{c} \in \{-1,1\}^n$ such that $c_1=1$, and let $A\in M_{n,n}(\mathbb{R})$.

Let us define a diagonal matrix $P(\vec{c})$ such that
\begin{equation}\label{eq:15}
P(\vec{c})=diag(c_1,c_2,\ldots,c_n)\text{.}
\end{equation}

The matrix $P(\vec{c})$ is well-known as a signature matrix.

Obviously,
\begin{equation}\label{eq:16}
P(\vec{c})^{-1}=P(\vec{c})\text{.}
\end{equation}

Furthermore, it is readily verified that
\begin{equation}\label{eq:17}
\varphi_{\vec{c}}(A)=P(\vec{c})\cdot A \cdot P(\vec{c})\text{.}
\end{equation}

By Eqns.~(\ref{eq:16}) and (\ref{eq:17}), it follows that a matrix $\varphi_{\vec{c}}(A)$ is similar to a matrix $A$.
Also, by the Eq.~(\ref{eq:17}), it follows that a matrix $\varphi_{\vec{c}}(A)$ is congruent to a matrix $A$. 

The Eq.~(\ref{eq:7.2}) is proved in Lemma \ref{l:2}.

This completes the proof of Theorem \ref{t:1}.

\section{Theorem \ref{t:1} in Terms of Operators}\label{sec:6}

Let $B=\{\vec{b_1},\vec{b_2},\ldots,\vec{b_n}\}$ be an arbitrary base of the vector space $\mathbb{R}^n$, and let $\varphi: \mathbb{R}^n \to \mathbb{R}^n$ be an operator such that has  matrix $A$ under the base $B$.

Let $\vec{c} \in \{-1,1\}^n$ . It is readily verified that the set $B(\vec{c})=\{c_1\vec{b_1},c_2\vec{b_2},\ldots,c_n\vec{b_n}\}$ is also a base of $\mathbb{R}^n$.

We assert that the matrix of an operator $\varphi$ under the new base $B(\vec{c})$ is matrix $\varphi_{\vec{c}}(A)$.

For $1\leq j \leq n$, we have that:

\begin{equation}\label{eq:18}
\varphi(\vec{b_j})=\sum_{i=1}^{n}a_{ij}\vec{b_i}\text{.}
\end{equation}

By using the Eq.~(\ref{eq:18}), we have gradually:
\begin{align}
\varphi(c_j\vec{b_j})&=c_j\cdot \varphi(\vec{b_j})\notag\\
&=\sum_{i=1}^{n}a_{ij}\cdot c_j\cdot \vec{b_i}\notag\\
&=\sum_{i=1}^{n}{c_i}^2 \cdot a_{ij}\cdot c_j \cdot \vec{b_i}\notag\\
&=\sum_{i=1}^{n}c_i \cdot a_{ij}\cdot c_j \cdot (c_i\vec{b_i})\label{eq:19}\text{.}
\end{align}

The Eq.~(\ref{eq:19}) proves our assertion.

\section{A Proof of Theorem \ref{t:2}}\label{sec:7}

Let $\vec{c}$ and $\vec{d} \in \{-1,1\}^n$ such that $c_1=d_1=1$.  Let us suppose that $\varphi_{\vec{c}}=\varphi_{\vec{d}}$.

Then, for any square matrix $A\in M_{n,n}(\mathbb{R})$,  $\varphi_{\vec{c}}(A)=\varphi_{\vec{d}}(A)$.

Let $B \in M_{n,n}(\mathbb{R})$ be a matrix such that $b_{1i}\neq 0$, for all $1 \leq i \leq n$. Then, it must be:

\begin{align}
(\varphi_{\vec{c}}(B))_{1i}&=(\varphi_{\vec{d}}(B))_{1i}\notag\\
c_1\cdot b_{1i}\cdot c_i&=d_1\cdot b_{1i}\cdot d_i\notag\\
b_{1i}\cdot c_i&= b_{1i}\cdot d_i \notag\\
c_i&=d_i\text{.}&&(\text{division by non-zero})\label{eq:20}
\end{align}

Hence, assumption $\varphi_{\vec{c}}=\varphi_{\vec{d}}$, yield to $\vec{c}=\vec{d}$.

Obviously, there are $2^{n-1}$ different maps $\varphi_{\vec{c}}$.

Furthermore, let us consider the composition of maps $\varphi_{\vec{d}} \circ \varphi_{\vec{c}}$.

We have gradually:

\begin{align}
(\varphi_{\vec{d}} \circ \varphi_{\vec{c}})(A)&=\varphi_{\vec{d}} (\varphi_{\vec{c}} (A))\notag\\
&=\varphi_{\vec{d}}(P(\vec{c})\cdot A\cdot  P(\vec{c}))&&(\text{by Eq.~(\ref{eq:17})})\notag\\
&=P(\vec{d})\cdot P(\vec{c})\cdot A\cdot  P(\vec{c})\cdot P(\vec{d})\label{eq:21}
\end{align}

It is well-known that any two diagonal matrix commute under the multiplication of matrices.

Let $\vec{e}\in \mathbb{R}^n$ be a vector such that $e_i=c_i\cdot d_i$, for all $1 \leq i \leq n$. Obviously, $\vec{e}\in \{-1,1\}^n$ and $e_1=1$.

It is readily verified that:

\begin{equation}
P(\vec{c})\cdot P(\vec{d})=P(\vec{d})\cdot P(\vec{c})=P(\vec{e})\text{.}\label{eq:22}
\end{equation}

By the Eq.~(\ref{eq:22}), the Eq.~(\ref{eq:21}) becomes
\begin{equation}\label{eq:23}
(\varphi_{\vec{d}} \circ \varphi_{\vec{c}})(A)=P(\vec{e})\cdot A \cdot P(\vec{e})\text{.}
\end{equation}

By Eqns.~(\ref{eq:17}) and (\ref{eq:23}), it follows that

\begin{equation}\label{eq:24}
\varphi_{\vec{c}} \circ \varphi_{\vec{d}}=\varphi_{\vec{d}} \circ \varphi_{\vec{c}}= \varphi_{\vec{e}}\text{.}
\end{equation}

Therefore, the set of all maps  $\Psi_n=\{\varphi_{\vec{c}}: \vec{c} \in \{-1,1\}^n \text{ and } c_1=1\}$ is closed under the composition of maps $\circ$.

Let $id: M_{n,n}(\mathbb{R})\to  M_{n,n}(\mathbb{R})$ denote identity map, i.e., $id(A)=A$, for any $A \in  M_{n,n}(\mathbb{R})$.

Let $f\in \{-1,1\}^n$ be a vector such that $f_i=1$ for all $1\leq i \leq n$. 

Obviously, $id=\varphi_{\vec{f}}$. Hence, $id \in \Psi_n$.

Furthermore, it is readily verified that

\begin{equation}
\varphi_{\vec{c}} \circ \varphi_{\vec{c}}=id\text{.}\label{eq:26}
\end{equation}

Therefore, the group $(\varPsi_n,\circ)$ is an abelian group, where every element has order two except of neutral element of the group.

Clearly, the group $(\varPsi_n,\circ) \cong ({\mathbb{Z}_2}^{n-1},+)$.

\section{An Example of The Group $(\varPsi_3,\circ)$}\label{sec:8}

Let $\vec{c}_1$, $\vec{c}_2$, $\vec{c}_3$, and $\vec{c}_4$ be vectors from the set $\{-1,1\}^3$ such that:
\begin{align*}
\vec{c}_1&=(1,1,-1)\\
\vec{c}_2&=(1,-1,1)\\
\vec{c}_3&=(1,-1,-1)\\
\vec{c}_4&=(1,1,1)\text{.}
\end{align*}

Let $A$ be an arbitrary matrix from the  $ M_{3,3}(\mathbb{R})$, i.e., 
\begin{equation}\notag
A=\begin{pmatrix}
a_{11}&a_{12}&a_{13}\\
a_{21}&a_{22}&a_{23}\\
a_{31}&a_{32}&a_{33}\text{.}
\end{pmatrix}
\end{equation}

Then $\varPsi_3=\{\varphi_{\vec{c_1}},\varphi_{\vec{c_2}},\varphi_{\vec{c_3}},\varphi_{\vec{c_4}}\}$, where these four maps are:

\begin{equation}\label{eq:26}
\varphi_{\vec{c_1}}(A)=\begin{pmatrix}
a_{11}&a_{12}&-a_{13}\\
a_{21}&a_{22}&-a_{23}\\
-a_{31}&-a_{32}&a_{33}\text{.}
\end{pmatrix}
\end{equation}

\begin{equation}\label{eq:27}
\varphi_{\vec{c_2}}(A)=\begin{pmatrix}
a_{11}&-a_{12}&a_{13}\\
-a_{21}&a_{22}&-a_{23}\\
a_{31}&-a_{32}&a_{33}\text{.}
\end{pmatrix}
\end{equation}

\begin{equation}\label{eq:28}
\varphi_{\vec{c_3}}(A)=\begin{pmatrix}
a_{11}&-a_{12}&-a_{13}\\
-a_{21}&a_{22}&a_{23}\\
-a_{31}&a_{32}&a_{33}\text{.}
\end{pmatrix}
\end{equation}

\begin{equation}\label{eq:29}
\varphi_{\vec{c_4}}(A)=\begin{pmatrix}
a_{11}&a_{12}&a_{13}\\
a_{21}&a_{22}&a_{23}\\
a_{31}&a_{32}&a_{33}\text{.}
\end{pmatrix}
\end{equation}

Obviously, 
$\varphi_{\vec{c_4}}=id$.

Cayley table for the group $\varPsi_3$ is given below:

\begin{center}
\begin{tabular}{ |c| c| c|c| c| }
\hline
$\circ$& $\varphi_{\vec{c_1}}$&$\varphi_{\vec{c_2}}$&$\varphi_{\vec{c_3}}$ & $id$\\
\hline
$\varphi_{\vec{c_1}}$ & $id$ & $\varphi_{\vec{c_3}}$ & $\varphi_{\vec{c_2}}$ & $\varphi_{\vec{c_1}}$\\
\hline
$\varphi_{\vec{c_2}}$& $\varphi_{\vec{c_3}}$ & $id$ &$\varphi_{\vec{c_1}}$&$\varphi_{\vec{c_2}}$\\
\hline
$\varphi_{\vec{c_3}}$& $\varphi_{\vec{c_2}}$ & $\varphi_{\vec{c_1}}$& $id$ &$\varphi_{\vec{c_3}}$\\
\hline
$id$& $\varphi_{\vec{c_1}}$ & $\varphi_{\vec{c_2}}$ & $\varphi_{\vec{c_3}}$ & $id$\\
\hline
\end{tabular}
\end{center}

Obviously, the group $(\varPsi_3,\circ)$ is isomorphic with the Klein-four group $(\mathbb{Z}_2\times \mathbb{Z}_2,+)$.

\section{Vector subspaces ($Sym_{\varphi_{\vec{c}}}$,+) and ($AntiSym_{\varphi_{\vec{c}}}$,+) for $n=3$ }\label{sec:8.1}

Let $\vec{c}_1$, $\vec{c}_2$, $\vec{c}_3$, and $\vec{c}_4$ be  the same vectors as in previous Section \ref{sec:8}. Let $B_i$ , and $C_i \in M_{3,3}(\mathbb{R})$, where $1\leq i \leq 4$.

It is readily verified that $B_1 \in Sym_{\varphi_{\vec{c_1}}}$ if and only if:

\begin{equation}\notag
B_1=
\begin{pmatrix}
a_{11}&a_{12}&0\\
a_{21}&a_{22}&0\\
0&0&a_{33}\text{.}
\end{pmatrix}
\end{equation}

Similarly, $C_1 \in AntiSym_{\varphi_{\vec{c_1}}}$ if and only if:

\begin{equation}\notag
C_1=
\begin{pmatrix}
0&0&a_{13}\\
0&0&a_{23}\\
a_{31}&a_{32}&0\text{.}
\end{pmatrix}
\end{equation}

Furthermore, $B_2 \in Sym_{\varphi_{\vec{c_2}}}$ if and only if:

\begin{equation}\notag
B_2=
\begin{pmatrix}
a_{11}&0&a_{13}\\
0&a_{22}&0\\
a_{31}&0&a_{33}\text{.}
\end{pmatrix}
\end{equation}

Similarly, $C_2 \in AntiSym_{\varphi_{\vec{c_2}}}$ if and only if:

\begin{equation}\notag
C_2=
\begin{pmatrix}
0&a_{12}&0\\
a_{21}&0&a_{23}\\
0&a_{32}&0\text{.}
\end{pmatrix}
\end{equation}

 A matrix $B_3 \in Sym_{\varphi_{\vec{c_3}}}$ if and only if:

\begin{equation}\notag
B_3=\begin{pmatrix}
a_{11}&0&0\\
0&a_{22}&a_{23}\\
0&a_{32}&a_{33}\text{.}
\end{pmatrix}
\end{equation}

 A matrix $C_3 \in Sym_{\varphi_{\vec{c_3}}}$ if and only if:

\begin{equation}\notag
C_3=\begin{pmatrix}
0&a_{12}&a_{13}\\
a_{21}&0&0\\
a_{31}&0&0\text{.}
\end{pmatrix}
\end{equation}

A matrix $B_4 \in Sym_{\varphi_{\vec{c_4}}}$ if and only if:

\begin{equation}\notag
B_3=\begin{pmatrix}
a_{11}&a_{12}&a_{13}\\
a_{21}&a_{22}&a_{23}\\
a_{31}&a_{32}&a_{33}\text{.}
\end{pmatrix}
\end{equation}

In other words, any  square matrix $ B \in M_{3,3}(\mathbb{R})$ is symmetric under the map $\varphi_{\vec{c_4}}$.

Finally, a matrix $C_4 \in Sym_{\varphi_{\vec{c_4}}}$ if and only if:

\begin{equation}\notag
C_4=\begin{pmatrix}
0&0&0\\
0&0&0\\
0&0&0\text{.}
\end{pmatrix}
\end{equation}

In other words, the only square matrix $ C\in M_{3,3}(\mathbb{R})$, antisymmetric under the map $\varphi_{\vec{c_4}}$, is zero matrix.

\begin{remark}\label{r:1.1}
Let $n$ be a fixed natural number greater than one.
Dimensions of subspaces $Sym_{\varphi_{\vec{c}}}$ and $AntiSym_{\varphi_{\vec{c}}}$ depend on number of ones of a vector $\vec{c}$.

Let $r$ denote the number of ones of a vector $\vec{c}$. Then:

\begin{align}
dim(Sym_{\varphi_{\vec{c}}})&=r^2+(n-r)^2\label{eq:35.1}\\
dim(AntiSym_{\varphi_{\vec{c}}})&=2r(n-r)\text{.}\label{eq:35.2}
\end{align}

\end{remark}

\section{A Proof of Theorem \ref{t:3}}\label{sec:9}

Let $\vec{c} \in \{-1,1\}^n$ such that $c_1=1$.

Let $A$ and $B$ be arbitrary matrices from the set $M_{n,n}(\mathbb{R})$, and let $\alpha$ be an arbitrary real number.

By the Eq.~(\ref{eq:1}), it is easy to show that a map $\varphi_{\vec{c}}$ is an operator, i.e., 

\begin{align}
\varphi_{\vec{c}}(A+B)=\varphi_{\vec{c}}(A)+\varphi_{\vec{c}}(B)\label{eq:30}\\
\varphi_{\vec{c}}(\alpha\cdot A)=\alpha\cdot \varphi_{\vec{c}}(A)\text{.}\label{eq:31}
\end{align}

Let $C, D \in Sym_{\vec{c}}$ , and let $\alpha, \beta$ be arbitrary real numbers.

We have gradually:

\begin{align}
\varphi_{\vec{c}}(\alpha C+\beta D)&=\varphi_{\vec{c}}(\alpha C)+\varphi_{\vec{c}}(\beta D)&&(\text{by the Eq.~(\ref{eq:30})})\notag\\
&=\alpha \varphi_{\vec{c}} (C)+\beta \varphi_{\vec{c}} (D)&&(\text{by the Eq.~(\ref{eq:31})})\notag\\
&=\alpha C+\beta D\text{.}&&(\text{by Def.~(\ref{def:1})})\label{eq:32}
\end{align}

The Eq.~(\ref{eq:32}) proves that $(Sym_{\vec{c}},+)$ is a vector subspace of the vector space $M_{n,n}(\mathbb{R})$ over the field of real numbers $\mathbb{R}$.

Similarly, it can be shown that $(AntiSym_{\vec{c}},+)$ is a vector subspace of the vector space $M_{n,n}(\mathbb{R})$ over the field of real numbers $\mathbb{R}$.

It is readily verified that  $A=A_S(\vec{c})+A_{AS}(\vec{c})$.

Let us prove that $A_S(\vec{c})\in Sym_{\vec{c}}$.

Let $i$ and $j$ be natural numbers such that $1\leq i,j\leq n$. 

By the Eq.~(\ref{eq:1}), it follows that

\begin{equation}\label{eq:33}
(\varphi_{\vec{c}}(A_S(\vec{c})))_{ij}=c_i\cdot (A_S(\vec{c}))_{ij}\cdot c_j\text{.}
\end{equation}

We have two cases!

The first case: Let us suppose that  $c_i\cdot c_j=1$.

Then the Eq.~(\ref{eq:33}) becomes:
\begin{equation}\label{eq:34}
(\varphi_{\vec{c}}(A_S(\vec{c})))_{ij}=A_S(\vec{c})_{ij}\text{.}
\end{equation}

The second case: Let us suppose that  $c_i\cdot c_j=-1$.

Then the Eq.~(\ref{eq:33}) becomes:
\begin{equation}\label{eq:35}
(\varphi_{\vec{c}}(A_S(\vec{c})))_{ij}=-A_S(\vec{c})_{ij}\text{.}
\end{equation}

By the Eq.~(\ref{eq:8}), it follows that $A_S(\vec{c})_{ij}=0$, for $c_i\cdot c_j=-1$.

Hence, $-A_S(\vec{c})_{ij}=A_S(\vec{c})_{ij}$. 

The Eq.~(\ref{eq:34}) follows from Eqns.~(\ref{eq:8}) and (\ref{eq:35}).

This proves that $A_S(\vec{c})\in Sym_{\vec{c}}$. 

Similarly, it can be shown that $A_{AS}(\vec{c})\in AntiSym_{\vec{c}}$.

Finally, let us prove that:
\begin{equation}\label{eq:36}
Sym_{\vec{c}} \cap AntiSym_{\vec{c}}=\{\bf{0}\}\text{,}
\end{equation}
where $ \bf{0}$ denotes a square matrix of an order $n$ which all entries are zeros.

Let $C\in M_{n,n}(\mathbb{R})$ such that $ C\in Sym_{\vec{c}} \cap AntiSym_{\vec{c}}$.

Then it must be simultaneously $\varphi_{\vec{c}}(C)=C$ and $\varphi_{\vec{c}}(C)=-C$.

It follows that $C=-C$, i. e., $C=\bf{0}$.

This proves the Eq.~(\ref{eq:36}) and  completes the proof of Theorem \ref{t:3}.

\begin{remark}\label{r:1.2}
Let $ A \in M_{n,n}(\mathbb{R})$.
Then $A_S(\vec{c})=\frac{1}{2} (A+\varphi_{\vec{c}}(A))$ and $A_{AS}(\vec{c})=\frac{1}{2} (A-\varphi_{\vec{c}}(A))$.
\end{remark}

\section{A Proof of Lemma \ref{l:3}}\label{sec:10}

Let us prove that

\begin{equation}\label{eq:37}
\varphi_{\vec{c}}(A\cdot B)=\varphi_{\vec{c}}(A)\cdot \varphi_{\vec{c}}(B)\text{,}
\end{equation}
where $A,B \in M_{n,n}(\mathbb{R})$.

By the Eq.~(\ref{eq:17}), it follows:
\begin{equation}\label{eq:38}
\varphi_{\vec{c}}(A\cdot B)=P(\vec{c})\cdot A \cdot B\cdot P(\vec{c})\text{.}
\end{equation}

By the Eq.~(\ref{eq:16}), it follows that Eq.~(\ref{eq:38}) becomes:
\begin{align*}
\varphi_{\vec{c}}(A\cdot B)&=(P(\vec{c})\cdot A \cdot P(\vec{c}))\cdot (P(\vec{c})\cdot B\cdot P(\vec{c}))\text{.}\\
&=\varphi_{\vec{c}}(A)\cdot \varphi_{\vec{c}}(B)\text{.}
\end{align*}

The last equation above proves the Eq.~(\ref{eq:37}).

Let us prove that $(Sym_{\vec{c}},\cdot)$ is closed under the multiplication of matrices.

Let $C, D \in Sym_{\vec{c}}$. By definition \ref{def:1}, we know that

\begin{align*}
\varphi_{\vec{c}}(C)&=C\\
\varphi_{\vec{c}}(D)&=D\text{.}
\end{align*}

By the Eq.~(\ref{eq:37}), it follows that

\begin{align}
\varphi_{\vec{c}}(C\cdot D)&=\varphi_{\vec{c}}(C)\cdot \varphi_{\vec{c}}(D)\notag\\
&=C \cdot D\text{.}\label{eq:39}
\end{align}

The Eq.~(\ref{eq:39}) proves that $(Sym_{\vec{c}},\cdot)$ is closed under the multiplication of matrices.

Let $I \in M_{n,n}(\mathbb{R})$ denote identity matrix, i. e., $I=diag\underset{\underbrace{n\text{- times}}}{(1,1,\ldots, 1)}$.

 Obviously, $I \in Sym_{\vec{c}}$. Moreover, $\varphi_{\vec{c}}(I)=I$, for all vectors $\vec{c}\in \{-1,1\}^n$ such that $c_1=1$.

By the Eq.~(\ref{eq:1}), it follows that if $A \in Sym_{\vec{c}}$ then  $-A \in Sym_{\vec{c}}$.

This completes the proof of Lemma \ref{l:3}.

\begin{remark}\label{r:1}
Let $\cdot$ denote standard multiplication of matrices.
Then the following rule is true:

\begin{center}
\begin{tabular}{ |c| c| c|}
\hline
$\cdot$ & $Sym_{\vec{c}}$ & $AntiSym_{\vec{c}}$\\
\hline
$Sym_{\vec{c}}$ & $Sym_{\vec{c}}$&$ AntiSym_{\vec{c}}$ \\
\hline
$AntiSym_{\vec{c}}$ & $AntiSym_{\vec{c}}$ & $Sym_{\vec{c}}$\\
\hline
\end{tabular}
\end{center}

This rule does not hold for multiplication of symmetric and antisymmetric matrices.

\end{remark}

\begin{remark}\label{r:2}

Let $\vec{c}\in \{-1,1\}^n$ such that $c_1=1$.

Let $Ker(\varphi_{\vec{c}})$ and $Im(\varphi_{\vec{c}})$ denote kernel and image of an operator $\varphi_{\vec{c}}$.
It is readily verified that  $Ker(\varphi_{\vec{c}})=\{\bf{0}\}$ and $Im(\varphi_{\vec{c}})=M_{n,n}(\mathbb{R})$. Hence, $\varphi_{\vec{c}}$ is a bijective operator.

By Eqns.~(\ref{eq:30}) and (\ref{eq:37}), it follows that $\varphi_{\vec{c}}$ is an automorphism of a ring $(M_{n,n}(\mathbb{R}),+,\cdot)$.
\end{remark}

\section{A Proof of Theorem \ref{t:4}}\label{sec:11}

Let $ A \in M_{n,n}(\mathbb{R})$, and let $ A(i,j)$ be an arbitrary principal minor of an order two of a matrix $A$ such that $1\leq i<j\leq n$.

By definition of a principal minor, it follows that

\begin{equation}\label{eq:40}
A(i,j)=\begin{vmatrix}
a_{ii}&a_{ij}\\
a_{ji}&a_{jj}
\end{vmatrix}
\end{equation}

Let us consider the following principal minors of an order two: $(A_S(\vec{c}))(i,j)$ and $(A_{AS}(\vec{c}))(i,j)$.
We assert that:
\begin{equation}\label{eq:41}
A(i,j)=A_S(\vec{c})(i,j)+A_{AS}(\vec{c})(i,j)\text{.}
\end{equation}

We have two cases.

The first case: let us suppose that $c_i\cdot c_j=1$.

By the Eq.~(\ref{eq:8}), it follows that

\begin{equation}\label{eq:42}
A_{S}(\vec{c})(i,j)=\begin{vmatrix}
a_{ii}&a_{ij}\\
a_{ji}&a_{jj}
\end{vmatrix}=A(i,j)
\end{equation}

By the Eq.~(\ref{eq:9}), it follows that

\begin{equation}\label{eq:43}
A_{AS}(\vec{c})(i,j)=\begin{vmatrix}
0&0\\
0&0
\end{vmatrix}=0
\end{equation}

By Eqns.~(\ref{eq:40}),(\ref{eq:42}), and (\ref{eq:43}), the Eq.~(\ref{eq:41}) follows.

The second case: let us suppose that $c_i\cdot c_j=-1$.

By the Eq.~(\ref{eq:8}), it follows that

\begin{equation}\label{eq:44}
A_{S}(\vec{c})(i,j)=\begin{vmatrix}
a_{ii}&0\\
0&a_{jj}
\end{vmatrix}=a_{ii}\cdot a_{jj}
\end{equation}

By the Eq.~(\ref{eq:9}), it follows that

\begin{equation}\label{eq:45}
A_{AS}(\vec{c})(i,j)=\begin{vmatrix}
0&a_{ij}\\
a_{ji}&0
\end{vmatrix}=-a_{ij}\cdot a_{ji}
\end{equation}

By Eqns.~(\ref{eq:40}), (\ref{eq:44}), and (\ref{eq:45}), the Eq.~(\ref{eq:41}) follows.

This proves the Eq.~(\ref{eq:41}).

By the Eq.~(\ref{eq:41}), it follows that:

\begin{align}
\underset{1\leq i <j \leq n}\sum A(i,j)&=\underset{1\leq i <j \leq n}\sum ( A_{S}(\vec{c})(i,j)+A_{AS}(\vec{c})(i,j))\notag\\
&=\underset{1\leq i <j \leq n}\sum A_{S}(\vec{c})(i,j)+\underset{1\leq i <j \leq n}\sum A_{AS}(\vec{c})(i,j)\text{.}\label{eq:46}
\end{align}

The Eq.~(\ref{eq:46}) proves Theorem \ref{t:4}.

By Eqns.~(\ref{eq:14}) and (\ref{eq:46}), the Eq.~(\ref{eq:11}) follows.

The proof of the Eq.~(\ref{eq:11.1}) is analoque to the proof of the Eq.~(\ref{eq:11}). Let $ A[i,j]$ be an arbitrary principal permanent of an order two of a matrix $A$ such that $1\leq i<j\leq n$.  It can be shown that:

\begin{equation}\label{eq:46.1}
A[i,j]=A_S(\vec{c})[i,j]+A_{AS}(\vec{c})[i,j]\text{.}
\end{equation}

By Eqns.~(\ref{eq:15.1}) and (\ref{eq:46.1}), the Eq.~(\ref{eq:11.1}) follows.

\section{A Proof of Theorem \ref{t:5}}\label{sec:12}

Let $ A \in M_{n,n}(\mathbb{R})$ and let $ A(i,j)$ be an arbitrary principal minor of an order two of a matrix $A$ such that $1\leq i<j\leq n$.

Let us consider the following principal minors of an order two: $A_S(i,j)$ and $A_{AS}(i,j)$.
We know that:

\begin{equation}\label{eq:47}
A_{S}(i,j)=\begin{vmatrix}
a_{ii}&\frac{1}{2}\cdot(a_{ij}+a_{ji})\\
\frac{1}{2}\cdot(a_{ij}+a_{ji})&a_{jj}
\end{vmatrix}
\end{equation}

\begin{equation}\label{eq:48}
A_{AS}(i,j)=\begin{vmatrix}
0&\frac{1}{2}\cdot(a_{ij}-a_{ji})\\
\frac{1}{2}\cdot(a_{ji}-a_{ij})&0
\end{vmatrix}
\end{equation}

By the Eq.~(\ref{eq:47}), it follows that:

\begin{equation}\label{eq:49}
A_{S}(i,j)=a_{ii}\cdot a_{jj}-\frac{1}{4}\cdot (a_{ij}+a_{ji})^2\text{.}
\end{equation}

By the Eq.~(\ref{eq:48}), it follows that:

\begin{equation}\label{eq:50}
A_{AS}(i,j)=\frac{1}{4}\cdot (a_{ij}-a_{ji})^2\text{.}
\end{equation}

By  adding Eqns.~(\ref{eq:49}) and (\ref{eq:50}), we have gradually:

\begin{align*}
A_{S}(i,j)+A_{AS}(i,j)&=a_{ii}\cdot a_{jj}-\frac{1}{4}\cdot (a_{ij}+a_{ji})^2+\frac{1}{4}\cdot (a_{ij}-a_{ji})^2\\
&=a_{ii}\cdot a_{jj}+\frac{1}{4}\cdot ((a_{ij}-a_{ji})^2-(a_{ij}+a_{ji})^2)\\
&=a_{ii}\cdot a_{jj}+\frac{1}{4}\cdot (a_{ij}-a_{ji}-a_{ij}-a_{ji}) \cdot (a_{ij}-a_{ji}+a_{ij}+a_{ji})\\
&=a_{ii}\cdot a_{jj}+\frac{1}{4}\cdot (-2a_{ji}) \cdot (2a_{ij})\\
&=a_{ii}\cdot a_{jj}- a_{ji} \cdot a_{ij}\\
&=A(i,j)\text{.}&&(\text{by the Eq.~(\ref{eq:40})})
\end{align*}

Therefore, we obtain that

\begin{equation}\label{eq:51}
A(i,j)=A_{S}(i,j)+A_{AS}(i,j)\text{.}
\end{equation}

By the Eq.~(\ref{eq:51}), it follows that

\begin{align}
\underset{1\leq i <j \leq n}\sum A(i,j)&=\underset{1\leq i <j \leq n}\sum ( A_{S}(i,j)+A_{AS}(i,j))\notag\\
&=\underset{1\leq i <j \leq n}\sum A_{S}(i,j)+\underset{1\leq i <j \leq n}\sum A_{AS}(i,j)\text{.}\label{eq:52}
\end{align}

The Eq.~(\ref{eq:52}) proves Theorem \ref{t:5}.

By Eqns.~(\ref{eq:14}) and (\ref{eq:52}), the Eq.~(\ref{eq:12}) follows.

The proof of the Eq.~(\ref{eq:12.1}) is analoque to the proof of the Eq.~(\ref{eq:12}). Let $ A[i,j]$ be an arbitrary principal permanent of an order two of a matrix $A$ such that $1\leq i<j\leq n$.  It can be shown that:

\begin{equation}\label{eq:52.1}
A[i,j]=A_S[i,j]+A_{AS}[i,j]\text{.}
\end{equation}

By Eqns.~(\ref{eq:15.1}) and (\ref{eq:52.1}), the Eq.~(\ref{eq:12.1}) follows.

\section{A Proof of Theorem \ref{t:6}}\label{sec:14}

Let $e_i \in M_{n,1}(\mathbb{R})$ form the standard base for the vector space $ M_{n,1}(\mathbb{R})$, where $1\leq i \leq n$.

Let us prove that vectors $A(e_l)$ and $e_s$ are orthogonal if $l \in I_r$ and $s\in J_{n-r}$.

Let us suppose that  $l \in I_r$ and $s\in J_{n-r}$. Then $c_l=1$ and $c_s=-1$. Hence, $c_l\cdot c_s=-1$.

The scalar product of vectors $A(e_l)$ and $e_s$ is equal to $a_{sl}$. Due to $c_l\cdot c_s=-1$ and $A \in Sym_{\vec{c}}$, it must be $a_{sl}=0$.

Therefore, vectors $A(e_l)$ and $e_s$ of  must be orthogonal if $l \in I_r$ and $s\in J_{n-r}$.

Similarly, vectors $e_l$ and $A(e_s)$ of  must be orthogonal if $l \in I_r$ and $s\in J_{n-r}$.

 Let $\varphi: M_{n,1}(\mathbb{R}) \to M_{n,1}(\mathbb{R})$ be an operator such that has  a matrix $A$ under the standard base $\{e_i: 1\leq i \leq n\}$ of the vector space $M_{n,1}(\mathbb{R})$.

Let $U$ denote the set $\{e_{i_1}, e_{i_2},\ldots e_{i_r}\}$ and let $V$  denote the set $\{e_{j_1}, e_{j_2},\ldots e_{j_{n-r}}\}$.

 Vector subspaces $span(U)$ and $span(V)$ are invariant under the operator $\varphi$.

Therefore, an operator $\varphi$ under the permutation of the standard base $\{e_{i_1}, e_{i_2},\ldots e_{i_r}, e_{j_1}, e_{j_2},\ldots e_{j_{n-r}}\}$ has the matrix $diag(D,E)$.

Hence, a matrix $A$ is similar to a block matrix $diag(D,E)$.

This proves Theorem \ref{t:6}.

Furthermore, it follows that:

\begin{equation}\label{eq:70.1}
P^{-1}\cdot A \cdot P=diag(D,E)\text{,}
\end{equation}
where $P$ is a permutation matrix.

For the proof of Corollary \ref{c:4}, see, for example, \cite{TzerHua}.

\section{A Proof of Theorem \ref{t:7}}\label{sec:15}

Let $e_i \in M_{n,1}(\mathbb{R})$ form the standard base for the vector space $ M_{n,1}(\mathbb{R})$, where $1\leq i \leq n$.

Let us prove that vectors $A(e_l)$ and $e_s$ are orthogonal if $l \in I_r$ and $s\in  I_r$.

Let us suppose that  $l \in I_r$ and $s\in I_r$. Then $c_l=1$ and $c_s=1$. Hence, $c_l\cdot c_s=1$.

The scalar product of vectors $A(e_l)$ and $e_s$ is equal to $a_{sl}$. Due to $c_l\cdot c_s=1$ and $A \in AntiSym_{\vec{c}}$, it must be $a_{sl}=0$.

Therefore, vectors $A(e_l)$ and $e_s$ of  must be orthogonal if $l \in I_r$ and $s\in I_r$.

Similarly, vectors $e_l$ and $A(e_s)$ of  must be orthogonal if $l \in J_{n-r}$ and $s\in J_{n-r}$.

Let $\varphi: M_{n,1}(\mathbb{R}) \to M_{n,1}(\mathbb{R})$ be an operator such that has  a matrix $A$ under the standard base $\{e_i: 1\leq i \leq n\}$ of the vector space $M_{n,1}(\mathbb{R})$.

Let $U$ denote the set $\{e_{i_1}, e_{i_2},\ldots e_{i_r}\}$ and let $V$  denote the set $\{e_{j_1}, e_{j_2},\ldots e_{j_{n-r}}\}$. 

Then $A(span (U))\subset span(V)$ and $A(span(V)) \subset span(U)$.

Therefore, an operator $\varphi$ under the same permutation of the  standard base:\\  $\{e_{i_1}, e_{i_2},\ldots e_{i_r}, e_{j_1}, e_{j_2},\ldots e_{j_{n-r}}\}$ has the matrix $H$ from the Eq.~(\ref{eq:16.4}).

This completes the proof of Theorem \ref{t:7}.

Furthermore, it follows that:

\begin{equation}\label{eq:70.2}
P^{-1}\cdot A \cdot P=H\text{,}
\end{equation}
where $P$ is the same permutation matrix from the Eq.~(\ref{eq:70.1}).

Corollary \ref{c:7} follows from Theorem \ref{t:7}. See , for example, \cite{TzerHua} for formulas about block matrices.

\section{A Proof of Theorem \ref{t:8}}\label{sec:15}

Let $A \in M_{n,n}(\mathbb{R})$ be a fixed matrix.

Let us consider the set $X_A=\{\varphi_{\vec{c}}(A):\vec{c}\in \{-1,1\}^n \text{ and } c_1=1\}$, and let $|X_A|$ denote the number of elements of the set $X_A$.

Let us suppose that  $\vec{c}\in \{-1,1\}^n$ be a fixed vector such that $c_1=1$. Let us consider the following set 
$ X_{A,\vec{c}}= \{\vec{d}\in \{-1,1\}^n : \varphi_{\vec{d}}(A)=\varphi_{\vec{c}}(A) \text{ and } d_1=1 \}$. We assert that $|X_{A,\vec{c}}|=2^{t-1}$.

Let $v_i$ and $v_j$ be two vertices of the graph $G$ which adjacency matrix is  matrix $A$ and $i <j$. We have two cases.

\textbf{The first case}: vertices $v_i$ and $v_j$ belong to the two different components of the graph $G$. Then, by Definition \ref{def:4}, $a_{ij}=a_{ji}=0$.
By the Eq.~(\ref{eq:1}), it follows that $(\varphi_{\vec{d}}(A))_{ij}=(\varphi_{\vec{c}}(A))_{ij}=0$. Obviously, the choice of $d_i$ and $d_j$ are independent from $c_i$ and $c_j$.

\textbf{The second case}: vertices $v_i$ and $v_j$ belong to the same component of the graph $G$.
Then there is a path $w$ between vertices $v_i$ and $v_j$. Let us suppose that a path $w$ is determined by vertices $v_{l_1}$,  $v_{l_2}$, $\ldots$, $v_{l_k}$, where $v_{l_1}=v_i$ and $v_{l_k}=v_j$ .
Without loss of generality, we may assume that all elements $a_{l_1 l_2}$, $a_{l_2 l_3}$, $\ldots$, $a_{l_{k-1} l_k}$ are non-zeros.

Let $s$ be a natural number such that $1 \leq s \leq k-1$. By using the Eq.~(\ref{eq:1}) and the fact $a_{l_s l_{s+1}}\neq 0$, from $(\varphi_{\vec{c}}(A))_{ l_s l_{s+1}}=(\varphi_{\vec{d}}(A))_{l_s  l_{s+1}}$, it follows that:

\begin{equation}\label{eq:70.3}
c_{l_s}\cdot c_{l_{s+1}}=d_{l_s}\cdot d_{l_{s+1}}\text{.}
\end{equation}

By using the Eq.~(\ref{eq:70.3}) , we have gradually:
\begin{align}
\prod_{s=1}^{k-1}c_{l_s}\cdot c_{l_{s+1}}&=\prod_{s=1}^{k-1}d_{l_s}\cdot d_{l_{s+1}}\notag\\
c_i\cdot \bigl{(}\prod _{s=1}^{k-1}{c_{l_s}}^2 \bigr{)}\cdot c_j&=d_i\cdot \bigl{(}\prod _{s=1}^{k-1}{d_{l_s}}^2\bigr{)}  \cdot d_j\notag\\
c_i\cdot c_j&=d_i\cdot d_j\text{.}\label{eq:70.4}
\end{align}

We have two cases.

The first case: $i\neq 1$.

We have two possibilities.

If $d_i=c_i$, then it must be $d_j=c_j$. Otherwise,  if $d_i=-c_i$, then it must be $d_j=-c_j$.

The second case: $i=1$.

Obviously,  $c_1=d_1=1$ and $ d_j=c_j$.

For each of $t-1$ connected components ( without vertex $v_1$) of a graph $G$, there are two choices for coordinates of a vector $\vec{d}$.

Hence, $|X_{A,\vec{c}}|=2^{t-1}$.

By using the  equation $ |X_A| \cdot 2^{t-1}=2^{n-1}$, it follows that:
\begin{equation}\label{eq:70.5}
|X_A|=2^{n-t}\text{.}
\end{equation}

The Eq.~(\ref{eq:70.5}) proves Theorem \ref{t:8}.

A proof of Corollary \ref{c:8} immediately follows from the proof of Theorem \ref{t:8}.

\section{Concluding Remarks}\label{sec:16}

\begin{remark}\label{r:5}
 Let us  define that a matrix $A \in M_{n,n}(\mathbb{R})$ is symmetric under the permutation matrix $P$ if  $P^{-1}\cdot A \cdot P=A$. Similarly,  let us  define that a matrix $A$ is antisymmetric  under the permutation matrix $P$ if  $P^{-1}\cdot A \cdot P=-A$.
Let $Sym_P$ denote the set of all symmetric matrix under the permutation matrix $P$ ,  and let $AntiSym_P$ denote the set of all antisymmetric matrix under the permutation matrix $P$. It is readily verified that $(Sym_P, +)$ and $(AntiSym_P,+)$ are vector subspaces of $( M_{n,n}(\mathbb{R}),+)$. However, the 
following equation 
\begin{equation}\label{eq:72}
Sym_P \oplus AntiSym_P=M_{n,n}(\mathbb{R})\text{,}
\end{equation}
is not true in general. 
\end{remark}

\section*{Acknowledgments}
I want to thank to professors Du\v{s}ko Bogdani\'{c}, Nik Stopar, and \DJ or\dj e Barali\'{c} for reading my manuscript. This research is not funded by any organization.

\end{document}